\documentclass[12pt]{article}
\title{Birational geometry of symplectic resolutions 
of nilpotent orbits II}
\author{Yoshinori Namikawa}
\date{ }
\usepackage{amsmath,amssymb,amsthm, amscd}
\chardef\bslash=`\\
\newtheorem{Thm}{Theorem}[section]

\newtheorem{Lem}[Thm]{Lemma}
\newtheorem{Prop}[Thm]{Proposition}
\newtheorem{Def}{Definition}
\newtheorem{Rque}[Thm]{Remark}

\newtheorem{Exam}[Thm]{Example}

\def\0{{\mathcal O}}
\def\g{{\mathfrak g}}

\def\h{{\mathfrak h}}

\begin{document}
\maketitle

\section{Introduction}
Let $G$ be a complex simple Lie group and let $\g$ be 
its Lie algebra. A nilpotent orbit $\mathcal{O}$ of $\g$ is an orbit 
of a nilpotent element $x \in \g$ by the adjoint action 
of $G$ on $\g$. The closure $\bar{\mathcal{O}}$ of $\mathcal{O}$ 
becomes a symplectic singularity via the Kostant-Kirillov form. 
By Fu [Fu], any symplectic 
resolution of $\bar{\mathcal{O}}$ is obtained as a Springer resolution 
$$ T^*(G/P) \to \bar{\mathcal{O}} $$ 
for a parabolic subgroup $P \subset G$.
In Part I [Na 2], when $\g$ is classical, 
we have proved that any two symplectic resolutions 
of $\bar{\mathcal{O}}$ 
are connected by a sequence of Mukai flops of type $A$ or of type $D$.  
In this paper (Part II), we shall improve and generalize all arguments in 
Part I so that the exceptional Lie algebras can be dealt with.  
We shall replace all arguments of [Na 2] which uses flags, by those 
which use only Dynkin diagrams.  
In the classical case, we already know which parabolic subgroups 
$P$ appear as the polarizations of $\mathcal{O}$ and when 
the Springer map $\mu : Y_P := T^*(G/P) \to \bar{\mathcal{O}}$ has 
degree 1 ([He]); so, in [Na 2], we only had to study the relationship 
between such polarizations. 
But, for the exceptional Lie algebras, no complete answer 
seems to be known. In this paper, we will start with a 
nilpotent orbit closure $\bar{\mathcal{O}}$ which has a 
Springer resolution $Y_{P_0} := T^*(G/P_0) \to \bar{\mathcal{O}}$.  
Even when $\g$ is classical, we will  
not use the classification of polarizations [He].   
First we introduce an equivalence relation in the set of 
parabolic subgroups of $G$ in terms of marked Dynkin diagrams 
(Definition 1). Our main theorem (Theorem 4.1) 
then claims that a parabolic subgroup 
$P \subset G$ always gives a Springer resolution of 
$\bar{\mathcal{O}}$ if $P$  
is equivalent to $P_0$. Moreover,  
any symplectic resolution of $\bar{\mathcal{O}}$ actually has 
this form, which will be proved as a corollary 
of the fact that the movable cone $\overline{\mathrm{Mov}}(Y_{P_0}/
\bar{\mathcal{O}})$ is the union of nef cones $\overline{\mathrm{Amp}}
(Y_P/\bar{\mathcal{O}})$ with $P \sim P_0$.  
Here all $Y_P$ ($P \sim P_0$) are connected 
by a sequence of certain Mukai flops (cf. Example 3.5, 
Theorem 4.1, (v)). When $\g$ is of type $E_6$, new Mukai flops 
(which are called of type $E_{6,I}$, $E_{6,II}$) appear.  
 
Finally, the author would like to thank S. Mukai for an 
important comment on [Na 2] and would like to thank D. Alvis 
for sending him the paper [Al].  

{\bf Notation}. For a proper birational map $f$ of algebraic 
varieties, we say that $f$ is {\em divisorial} if 
$\mathrm{Exc}(f)$ contains a divisor, and otherwise, we 
say that $f$ is {\em small}. Note that the terminology of 
"small" is, for example, different from that in [B-M].  

\section{Nilpotent orbits and Springer's correspondence}
Let $G$ be a complex simple Lie group and let $B$ be a Borel 
subgroup of $G$. Let $\g$ (resp. $\mathfrak{b}$) 
be the Lie algebra of $G$ (resp. $B$). 
The set of nilpotent elements $\mathcal{N}$ of $\g$ 
is called the nilpotent variety. It coincides with the 
closure of the regular nilpotent orbit in $\g$. 
The (original) Springer resolution 
$$ \pi: T^*(G/B) \to \mathcal{N} $$ 
is constructed as follows. Let $n(\mathfrak{b})$ be 
the nil-radical of $\mathfrak{b}$. Then the cotangent 
bundle $T^*(G/B)$ of $G/B$ is identified with 
$G \times^B n(\mathfrak{b})$, which is, by definition, 
the quotient space of $G \times n(\mathfrak{b})$ by 
the equivalence relation $\sim$. Here $(g,x) \sim 
(g', x')$ if $g = gb$ and $x' = Ad_{b^{-1}}(x)$ for 
some $b \in B$. Then we define 
$\pi([g,x]) := Ad_g(x).$ 
According to Borho-MacPherson [B-M], we shall briefly 
review Springer's correspondence [Sp]. The nilpotent 
variety $\mathcal{N}$ is decomposed into the 
disjoint union of nilpotent orbits $\mathcal{O}_x$, where 
$x$ is a distinguished base point of the orbit 
$\mathcal{O}_x$. We put $d_x := \dim \pi^{-1}(x)$. 
Now $\pi_1(\mathcal{O}_x)$ acts on $H^{2d_x}(\pi^{-1}(x), 
\mathbf{Q})$ by monodromy. Decompose $H^{2d_x}(\pi^{-1}(x), 
\mathbf{Q})$ into irreducible representations of 
$\pi_1(\mathcal{O}_x)$:  
$$ H^{2d_x}(\pi^{-1}(x), \mathbf{Q}) = \oplus_{\phi}
(V_{\phi}\otimes V_{(x,\phi)}), $$ 
where $\phi : \pi_1(\mathcal{O}_x) \to 
\mathrm{End}(V_{\phi})$ are irreducible representations and 
$V_{(x,\phi)} = \mathrm{Hom}_{\pi_1(\mathcal{O}_x)}
(V_{\phi}, H^{2d_x}(\pi^{-1}(x), \mathbf{Q}))$. By 
definition, $\dim V_{(x, \phi)}$ coincides with the 
multiplicity of $\phi$ in $H^{2d_x}(\pi^{-1}(x), \mathbf{Q})$. 
We call $(x,\phi)$ is $\pi$-{\em relevant} if 
$V_{(x,\phi)} \ne \emptyset.$      
Fix a maximal torus $T$ in $B$, and let $W$ be the Weyl 
group relative to $T$. Then there is a natural action 
of $W$ on $H^{2d_x}(\pi^{-1}(x), \mathbf{Q})$ commuting with 
the action of $\pi_1(\mathcal{O}_x)$. Each factor 
$V_{\phi}\otimes V_{(x,\phi)}$ becomes a $W$-module, where 
$W$ acts trivially on $V_{\phi}$ and $V_{(x,\phi)}$ is 
an irreducible representation of $W$. These representations 
were originally constructed by Springer. In [B-M], they 
are given in terms of the decomposition 
theorem of intersection cohomology by Beilinson, Bernstein, 
Deligne and Gabber. The following theorem is called 
Springer's correspondence: 

\begin{Thm}
Any irreducible representaion of $W$ is isomorphic 
to $V_{(x,\phi)}$ for a unique $\pi$-relevant pair 
$(x,\phi)$. 
\end{Thm}

One can find the tables on Springer's correspondence 
in [C, 13.3] for each simple Lie group (see also 
[A-L], [B-L]).

\section{Parabolic subgroups and marked Dynkin diagrams}
Let $G$ be a complex reductive Lie group and let $\g$ 
be its Lie algebra. Fix a Cartan subalgebra $\h$ of $\g$ and let 
$$\g = \h \oplus \bigoplus_{\alpha \in \Phi}\g_{\alpha}$$ 
be the root space decomposition. Let $\Delta \subset \Phi$ 
be a base of $\Phi$ and denote by $\Phi^+$ (resp. $\Phi^-$) 
the set of positive roots (resp. negative root).   
We define a Borel subalgebra 
$\mathfrak{b}$ of $\g$ as    
$$ \mathfrak{b} := 
\h \oplus \bigoplus_{\alpha \in \Phi^+}\g_{\alpha}.$$ 
For a subset $\Theta \subset \Delta$, let $<\Theta>$ be 
the sub-root system generated by $\Theta$. We put
$<\Theta>^+ := <\Theta> \cap \Phi^+$ and 
$<\Theta>^- := <\Theta> \cap \Phi^-$. 
We define    
$$ \mathfrak{p}_{\Theta} := 
\mathfrak{h} \oplus \bigoplus_{\alpha \in \Phi^+}\g_{\alpha} 
\oplus \bigoplus_{\alpha \in <\Theta>^-}\g_{\alpha}. $$ 
By definition, $\mathfrak{p}_{\Theta}$ is a parabolic 
subalgebra containing $\mathfrak{b}$. Moreover, any parabolic 
subalgebra $\mathfrak{p}$ of $\g$ is $G$-conjugate to 
$\mathfrak{p}_{\Theta}$ for some $\Theta \subset \Delta$. 
$\mathfrak{p}_{\Theta}$ and $\mathfrak{p}_{\Theta'}$ are 
$G$-conjugate if and only if $\Theta = \Theta'$.  
Therefore, there is a one-to-one correspondence between 
subsets of $\Delta$ and the conjugacy classes of parabolic 
subalgebras of $\g$.  
An element of $\Delta$ is called a simple root, which 
corresponds to a vertex of the Dynkin diagram attached to 
$\g$. A Dynkin diagram with some vertices being marked is 
called a {\em marked Dynkin diagram}. If $\Theta \subset 
\Delta$ is given, we have a marked Dynkin diagram by marking 
the vertices which correspond to $\Delta \setminus \Theta$.  
A marked Dynkin diagram with only one marked vertex is called 
a {\em single marked Dynkin diagram}. A conjugacy class of 
parabolic subgroups $P \subset G$ with $b_2(G/P) = 1$ 
corresponds to a single marked Dynkin diagram. 

\begin{Exam}
When $G = SL(n)$, the parabolic subgroup of flag type 
$(k,n-k)$ corresponds  
to the marked Dynkin diagram 

\begin{picture}(300,20)(0,0) 
\put(30,0){\circle{5}}\put(35,0){\line(1,0){25}} 
\put(65,-3.5){- - -}\put(90,0){\line(1,0){15}}
\put(105,0){\circle*{5}}\put(110,0){\line(1,0){10}}
\put(100,-10){k}\put(125,-3.5){- - -}\put(150,0)
{\line(1,0){55}}\put(210,0){\circle{5}.}    
\end{picture} 
\vspace{0.4cm}
 
When $G = SO(2n+1)$, the parabolic subgroup of flag type 
$(k,2n-2k+1,k)$ corresponds  
to the marked Dynkin diagram 

\begin{picture}(300,20)(0,0) 
\put(30,0){\circle{5}}\put(35,0){\line(1,0){25}} 
\put(65,-3.5){- - -}\put(90,0){\line(1,0){15}}
\put(105,0){\circle*{5}}\put(110,0){\line(1,0){10}}
\put(100,-10){k}\put(125,-3.5){- - -}\put(150,0)
{\line(1,0){15}}\put(170,0){\circle{5}}\put(175,-3){$\Rightarrow$}
\put(190,0){\circle{5}.}    
\end{picture}
\vspace{0.4cm}
 
When $G = Sp(2n)$, the parabolic subgroup of 
flag type $(k,2n-2k,k)$ corresponds to  
the marked Dynkin diagram 

\begin{picture}(300,20)(0,0) 
\put(30,0){\circle{5}}\put(35,0){\line(1,0){25}} 
\put(65,-3.5){- - -}\put(90,0){\line(1,0){15}}
\put(105,0){\circle*{5}}\put(110,0){\line(1,0){10}}
\put(100,-10){k}\put(125,-3.5){- - -}\put(150,0)
{\line(1,0){15}}\put(170,0){\circle{5}}\put(175,-3){$\Leftarrow$}
\put(190,0){\circle{5}.}    
\end{picture}
\vspace{0.4cm}
 
Finally, assume that $G = SO(2n)$. Then the parabolic subgroup 
corresponding to the marked Dynkin diagram ($k \geq 3$) 

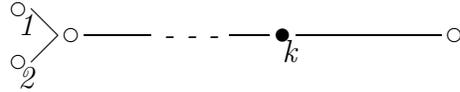
\begin{picture}(300,20)(0,0) 
\put(30,10){\circle{5}}\put(30,0){1}\put(30,-10){\circle{5}}
\put(30,-20){2}
\put(35,10){\line(1,-1){10}}\put(35,-10){\line(1,1){10}}
\put(50,0){\circle{5}}\put(55,0){\line(1,0){25}}
\put(85,-3.5){- - -}\put(110,0){\line(1,0){15}}\put(130,0)
{\circle*{5}}\put(130,-10){k}\put(135,0){\line(1,0){55}}
\put(195,0){\circle{5}}
\end{picture}
\vspace{0.7cm} 

has flag type $(n-k+1, 2k-2, n-k+1)$. On the other hand, 
two marked Dynkin 
diagrams 

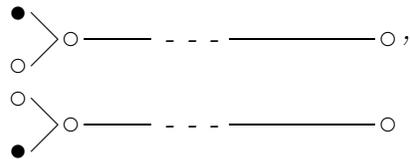
\begin{picture}(300,20)(0,0) 
\put(30,10){\circle*{5}}\put(30,-10){\circle{5}}
\put(35,10){\line(1,-1){10}}\put(35,-10){\line(1,1){10}}
\put(50,0){\circle{5}}\put(55,0){\line(1,0){25}}
\put(85,-3.5){- - -}\put(110,0){\line(1,0){55}}\put(170,0)
{\circle{5},}
\end{picture} 
\vspace{0.4cm}

\begin{picture}(300,20)(0,0) 
\put(30,10){\circle{5}}\put(30,-10){\circle*{5}}
\put(35,10){\line(1,-1){10}}\put(35,-10){\line(1,1){10}}
\put(50,0){\circle{5}}\put(55,0){\line(1,0){25}}
\put(85,-3.5){- - -}\put(110,0){\line(1,0){55}}\put(170,0)
{\circle{5}}
\end{picture}       
\vspace{0.4cm}

both give parabolic subgroups of flag type $(n,0,n)$
which are not $G$-conjugate. 
\end{Exam}
\vspace{0.2cm}

For a parabolic subgroup $P$ of $G$, let $\mathfrak{p}$ be 
its Lie algebra and let $n(\mathfrak{p})$ be the nil-radical 
of $\mathfrak{p}$. There is a unique nilpotent orbit $\mathcal{O} 
\subset \g$ such that $\mathcal{O} \cap n(\mathfrak{p})$ is 
an open dense subset of $n(\mathfrak{p})$. This nilpotent orbit 
is called the {\em Richardson orbit} for $P$. The cotangent 
bundle $T^*(G/P)$ of the homogenous space $G/P$ is naturally 
isomorphic to $G \times^P n(\mathfrak{p})$, which is the quotient 
space of $G \times n(\mathfrak{p})$ by the equivalence relation 
$\sim$. Here $(g,x) \sim (g', x')$ if $g' = gp$ and $x' = Ad_{p^{-1}}
(x)$ for some $p \in P$.  The {\em Springer map}  
$$ \mu: T^*(G/P) \to \overline{\mathcal{O}} $$ 
is defined as $\mu ([g,x]) = Ad_g(x)$. The Springer map $\mu$ 
is a generically finite surjective proper map. When 
$\mathrm{deg}\mu = 1$, it is called a Springer resolution. 
For a nilpotent orbit $\mathcal{O}_x \subset \overline{\mathcal{O}}$, 
we call $\mathcal{O}_x$ is $\mu$-relevant if 
$$\dim \mu^{-1}(x) = \mathrm{codim}(\mathcal{O}_x \subset 
\overline{\mathcal{O}})/2.$$ 
For the Springer resolution $\pi$ for a Borel subgroup $B$, 
every nilpotent orbit is $\pi$-relevant. However, this is not 
the case for a general parabolic subgroup $P$. The $\mu$-relevancy 
is closely related to Springer's correspondence. In order to 
state the result, we shall prepare some terminology. 
Let $L$ be a Levi subgroup of $P$. Fix a maximal torus $T$ of $L$. 
Then $T$ is also a maximal torus of $G$. Let $W(L)$ be the 
Weyl group for $L$ relative to $T$ and let $W$ be the Weyl group 
for $G$ relative to $T$. Now we have a natural inclusion 
$W(L) \subset W$. Let $\epsilon_{W(L)}$ be the sign representation 
of $W(L)$. Denote by $\epsilon_{W(L)}^W$ the induced representation 
of $\epsilon_{W(L)}$ to $W$. By section 1, every irreducible 
representation of $W$ has the form $V_{(x,\phi)}$ for a $\pi$-relevant 
pair $(x, \phi)$. Recall that $\phi$ is an irreducible representation 
of $\pi_1(\mathcal{O}_x)$. Denote by $1$ the trivial representation. 
Then $(x,1)$ is a $\pi$-relevant pair (cf. [B-M, Lemma 1.2]).  
  
\begin{Prop}
A nilpotent orbit $\mathcal{O}_x \subset \overline{\mathcal{O}}$ 
is $\mu$-relevant if and only if $V_{(x,1)}$ occurs 
in $\epsilon_{W(L)}^W$. 
\end{Prop}

{\em Proof}. See [B-M, Collorary 3.5, (b)].  
\vspace{0.15cm}

\begin{Prop}\label{cont} 
Let $G$ be a complex simple Lie group.  
Assume that $b_2(G/P) = 1$. Then the following are equivalent. 
\vspace{0.12cm}

(i) $\mathrm{deg}\mu = 1$ and $\mathrm{Codim}(\mathrm{Exc}(\mu)) 
\geq 2$,  \vspace{0.12cm} 

(ii) The single marked Dynkin diagram associated with $P$ is one 
of the following: \vspace{0.15cm}  

$A_{n-1}$ $(k < n/2)$ 
  
\begin{picture}(300,20)(0,0) 
\put(30,0){\circle{5}}\put(35,0){\line(1,0){25}} 
\put(65,-3.5){- - -}\put(90,0){\line(1,0){15}}
\put(105,0){\circle*{5}}\put(110,0){\line(1,0){10}}
\put(100,-10){k}\put(125,-3.5){- - -}\put(150,0)
{\line(1,0){55}}\put(210,0){\circle{5}}    
\end{picture}  

\begin{picture}(300,20)(0,0) 
\put(30,0){\circle{5}}\put(35,0){\line(1,0){25}} 
\put(65,-3.5){- - -}\put(90,0){\line(1,0){15}}
\put(105,0){\circle*{5}}\put(110,0){\line(1,0){10}}
\put(100,-10){n-k}\put(125,-3.5){- - -}\put(150,0)
{\line(1,0){55}}\put(210,0){\circle{5}}    
\end{picture}
\vspace{0.15cm}
 
$D_n$ $(n:$ $\mathrm{odd} \geq 4)$  

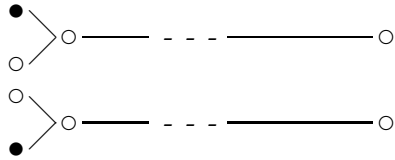
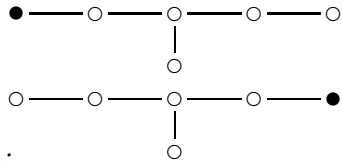
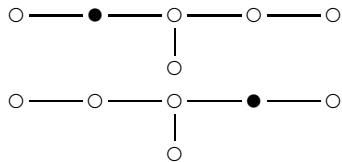
\begin{picture}(300,20)(0,0) 
\put(30,10){\circle*{5}}\put(30,-10){\circle{5}}
\put(35,10){\line(1,-1){10}}\put(35,-10){\line(1,1){10}}
\put(50,0){\circle{5}}\put(55,0){\line(1,0){25}}
\put(85,-3.5){- - -}\put(110,0){\line(1,0){55}}\put(170,0)
{\circle{5}}
\end{picture} 
\vspace{0.4cm}

\begin{picture}(300,20)(0,0) 
\put(30,10){\circle{5}}\put(30,-10){\circle*{5}}
\put(35,10){\line(1,-1){10}}\put(35,-10){\line(1,1){10}}
\put(50,0){\circle{5}}\put(55,0){\line(1,0){25}}
\put(85,-3.5){- - -}\put(110,0){\line(1,0){55}}\put(170,0)
{\circle{5}}
\end{picture}       
\vspace{0.4cm}

$E_{6,I}$:  

\begin{picture}(300,20)
\put(30,0){\circle*{5}}\put(35,0){\line(1,0){20}}
\put(60,0){\circle{5}}\put(65,0){\line(1,0){20}}
\put(90,0){\circle{5}}\put(90,-5){\line(0,-1){10}}
\put(90,-20){\circle{5}}\put(95,0){\line(1,0){20}}
\put(120,0){\circle{5}}\put(125,0){\line(1,0){20}}
\put(150,0){\circle{5}} 
\end{picture} 
\vspace{0.4cm}

\begin{picture}(300,20)
\put(30,0){\circle{5}}\put(35,0){\line(1,0){20}}
\put(60,0){\circle{5}}\put(65,0){\line(1,0){20}}
\put(90,0){\circle{5}}\put(90,-5){\line(0,-1){10}}
\put(90,-20){\circle{5}}\put(95,0){\line(1,0){20}}
\put(120,0){\circle{5}}\put(125,0){\line(1,0){20}}
\put(150,0){\circle*{5}} 
\end{picture}
\vspace{0.4cm}

$E_{6,II}$: 

\begin{picture}(300,20)
\put(30,0){\circle{5}}\put(35,0){\line(1,0){20}}
\put(60,0){\circle*{5}}\put(65,0){\line(1,0){20}}
\put(90,0){\circle{5}}\put(90,-5){\line(0,-1){10}}
\put(90,-20){\circle{5}}\put(95,0){\line(1,0){20}}
\put(120,0){\circle{5}}\put(125,0){\line(1,0){20}}
\put(150,0){\circle{5}} 
\end{picture}
\vspace{0.4cm}

\begin{picture}(300,20)
\put(30,0){\circle{5}}\put(35,0){\line(1,0){20}}
\put(60,0){\circle{5}}\put(65,0){\line(1,0){20}}
\put(90,0){\circle{5}}\put(90,-5){\line(0,-1){10}}
\put(90,-20){\circle{5}}\put(95,0){\line(1,0){20}}
\put(120,0){\circle*{5}}\put(125,0){\line(1,0){20}}
\put(150,0){\circle{5}} 
\end{picture}
\end{Prop}  
\vspace{0.4cm}

\begin{Rque}
In (ii) there are exactly two different markings for 
each Dynkin diagram $A_{n-1}$ with $k < n/2$, 
$D_n$, $E_{6,I}$ or $E_{6,II}$. They are called {\em dual} marked 
Dynkin diagrams. Let $P$ and $P'$ be the corresponding 
(conjugacy classes of) parabolic subgroups of $G$. Then $\mathfrak{p}$ 
and $\mathfrak{p}'$ have conjugate Levi factors. This 
implies that $P$ and $P'$ have the same Richardson orbit.  
\end{Rque}  

{\em Proof of Proposition 3.3}. Assume that the single marked Dynkin diagram 
is one of first two series in (ii). Then, by 
[Na 2], Lemmas 3.1 and 3.3, we already know 
that the Springer map $\mu: T^*(G/P) \to \bar{\mathcal{O}}$ 
becomes a small resolution (cf. Notation). 
If the single marked diagram 
is of type $E_{6,I}$, then the Richardson orbit $\mathcal{O}$ 
of $P$ coincides with orbit $\mathcal{O}_{2A_1}$ 
in the list of [C-M],p.129, which has dimension 
32. The maximal orbit contained in $\bar{\mathcal{O}}_{2A_1}
- \mathcal{O}_{2A_1}$ is $\mathcal{O}_{A_1}$, which has 
dimension 22. This shows that $\mathrm{Sing}(\bar{\mathcal{O}})$ 
has codimension $\geq$ 10 in $\bar{\mathcal{O}}$. 
On the other hand, since $\pi_1(\mathcal{O}_{2A_1}) = 1$ (cf. [C-M], 
p.129), $\mathrm{deg}(\mu) = 1$. If $\mu$ is a divisorial 
birational contraction, then $\mathrm{Codim}(\mathrm{Sing}
(\bar{\mathcal{O}} \subset \bar{\mathcal{O}}) = 2$ (cf. 
[Na 1, Cor. 1.5]), which is 
absurd. Hence $\mu$ should be a small resolution. 
If the single marked diagram is of type $E_{6,II}$, then 
the Richardson orbit $\mathcal{O}$ of $P$ coincides with 
the orbit $\mathcal{O}_{A_2 + 2A_1}$ in the list of 
[C-M], p.129, which has dimension 50. 
Moreover, $\pi_1(\mathcal{O}_{A_2 + 2A_1}) = 1$. 
By looking at the closure ordering of $E_6$ orbits 
(cf. [C], p.441), we see that the maximal orbit contained 
in $\bar{\mathcal{O}}_{A_2 + 2A_1} - \mathcal{O}_{A_2 + 2A_1}$ 
is the orbit $\mathcal{O}_{A_2 + A_1}$, which has dimension 
46. By the same argument as above, $\mu$ becomes a small 
resolution.       

To prove the implication (i) $\Rightarrow$ (ii), 
let us assume that the 
single marked Dynkin diagram is not contained in the list 
of (ii). Let $\mathcal{O}$ be the corresponding 
Richardson orbit. We shall first prove that $\bar{\mathcal{O}}$ 
contains a nilpotent orbit $\mathcal{O}'$ of 
codimension 2 (STEP 1). Next we shall prove that $\mathcal{O}'$ 
is $\mu$-relevant(STEP 2). These imply that $\mu$ is a 
divisorial birational contraction map if $\mathrm{deg}(\mu) = 1$. 

STEP 1: Assume that $\g$ is classical. If $\g$ is 
of type $A_{n-1}$, then we must look at the single 
marked Dynkin diagram with $k = n/2$. In this case, 
we already know that $\mu$ is a divisorial birational 
contraction map by [Na 2, Remark 3.2].  

When $\g$ is of type $B_n$, $C_n$ or $D_n$, the parabolic subgroup  
$P$ is a stabilizer group of an admissible isotropic flag. 
Its flag type is written as $(k,q,k)$. When $\g$ is of type 
$B_n$, we have $k > 0$, $q > 0$ and $2k + q = 2n + 1$. 
When $\g$ is of type $C_n$ or of type $D_n$, we have 
$k > 0$, $q \geq 0$ and $2k + q = 2n$.    
Denote by $\pi$ the dual partition of $\mathrm{ord}
(k,q,k)$ and call $\pi$ the Levi type of $P$. 

Assume that $\g$ is of type $B_n$. The Levi type of $P$  
is given by 

$$\pi := \left\{ 
\begin{array}{rl}    
\mathrm{[}3^{2n+1-2k}, 2^{3k-2n-1}\mathrm{]} & \quad (k > 1/3(2n+1))\\ 
\mathrm{[}3^k, 1^{2n-3k+1}\mathrm{]}  & \quad (k \leq 1/3(2n+1)) 
\end{array}\right.$$    

When $k > 1/3(2n+1)$, $k$ must be an odd number. In fact, 
if $k$ is even, then $I(\pi) \ne \emptyset$ and $\mathrm{deg}
(\mu) > 1$ (cf. [Na 2, Theorem 2.8]).  
Recall that the Richardson orbit $\mathcal{O}$ of $P$ has the Jordan type 
$S(\pi)$, where $S$ is the Spaltenstein map (cf. [Na 2, 
Theorem 2.7]). Since now $I(\pi) = \emptyset$, $S(\pi) 
= \pi$. Let us consider the nilpotent orbit $\mathcal{O}'$ 
of the Jordan type $[3^{2n+1-2k}, 2^{3k-2n-3}, 1^4]$ 
(resp. $[3^{k-1}, 2^2, 1^{2n-3k}]$, $[3^{k-1}, 1^3]$) 
when $k > 1/3(2n+1)$ (resp. $k < 1/3(2n+1)$, $k = 1/3(2n+1)$). 
In any case, we have $\mathcal{O}' \subset \bar{\mathcal{O}}$. 
By the dimension formula of nilpotent orbits ([C-M, Corollary 
6.1.4]), we see that $\dim \mathcal{O}' = \dim \mathcal{O} - 2$. 
 
Assume that $\g$ is of type $C_n$. The Levi type of $P$ 
is given by 

$$\pi := \left\{ 
\begin{array}{rl}    
\mathrm{[}3^{2n-2k}, 2^{3k-2n}\mathrm{]} & \quad (k > 2n/3)\\ 
\mathrm{[}3^k, 1^{2n-3k}\mathrm{]}  & \quad (k \leq 2n/3) 
\end{array}\right.$$
 
When $k \leq 2n/3$, $k$ must be an even number. In fact, 
if $k$ is odd, then $I(\pi) \ne \emptyset$ and $\mathrm{deg}
(\mu) > 1$ (cf. [Na 2, Theorem 2.8]). 
The Richardson orbit $\mathcal{O}$ has the Jordan type $\pi$. 
Let us consider the nilpotent orbit $\mathcal{O}'$ of 
the Jordan type $[3^{2n-2k}, 2^{3k-2n-1}, 1^2]$ (resp. 
$[3^{k-2}, 2^4, 1^{2n-3k-2}]$, $[3^{k-2}, 2^3]$) 
when $k > 2n/3$ (resp. $k < 2n/3$, $k = 2n/3$).  
In any case, we have $\mathcal{O}' \subset \bar{\mathcal{O}}$. 
By the dimension formula of nilpotent orbits ([C-M, Corollary 
6.1.4]), we see that $\dim \mathcal{O}' = \dim \mathcal{O} - 2$.

Assume that $\g$ is of type $D_n$. First assume that the Levi type of $P$  
is $[2^k]$. The single marked Dynkin diagram is not contained in the 
list of (ii) exactly when $k$ is even. In this case, $\mu$ is 
a divisorial birational contraction map by [Na 2, Remark 3.4]. 
We next assume $k < n$. In this case, the Levi type of $P$ is 
given by  

$$\pi := \left\{ 
\begin{array}{rl}    
\mathrm{[}3^{2n-2k}, 2^{3k-2n}\mathrm{]} & \quad (n > k > 2n/3)\\ 
\mathrm{[}3^k, 1^{2n-3k}\mathrm{]}  & \quad (k \leq 2n/3) 
\end{array}\right.$$ 

When $k > 2n/3$, $k$ must be an even number. In fact, 
if $k$ is odd, then $I(\pi) \ne \emptyset$ and $\mathrm{deg}
(\mu) > 1$ (cf. [Na 2, Theorem 2.8]).  
Recall that the Richardson orbit $\mathcal{O}$ of $P$ has the Jordan type 
$S(\pi)$, where $S$ is the Spaltenstein map (cf. [Na 2, 
Theorem 2.7]). Since now $I(\pi) = \emptyset$, $S(\pi) 
= \pi$. Let us consider the nilpotent orbit $\mathcal{O}'$ 
of the Jordan type $[3^{2n-2k}, 2^{3k-2n-2}, 1^4]$ 
(resp. $[3^{k-1}, 2^2, 1^{2n-3k-1}]$, $[3^{k-1}, 1^3]$) 
when $k > 2n/3$ (resp. $k < 2n/3$, $k = 2n/3$). 
In any case, we have $\mathcal{O}' \subset \bar{\mathcal{O}}$. 
By the dimension formula of nilpotent orbits ([C-M, Corollary 
6.1.4]), we see that $\dim \mathcal{O}' = \dim \mathcal{O} - 2$. 
 
When $\g$ is of type $G_2$, there are exactly two single marked 
Dynkin diagrams. In the table of $G_2$ nilpotent orbits in 
[C-M, p.128], $\mathcal{O}_{G_2(a_1)}$ is the Richardson orbit 
of the parabolic subgroups corresonding to these diagrams. 
The orbit $\mathcal{O}_{\tilde{A}_1}$ is contained in $\bar{\mathcal{O}}
_{G_2(a_1)}$. Note that $\dim \mathcal{O}_{G_2(a_1)} = 10$ and 
$\dim \mathcal{O}_{\tilde{A}_1} = 8$.  

When $\g$ is of type $F_4$, there are exactly four single marked 
Dynkin diagrams. Richardson orbits of the parabolic subgroups 
corresponding to them are $\mathcal{O}_{A_2}$, $\mathcal{O}_{\tilde{A}_2}$, 
$\mathcal{O}_{F_4(a_3)}$ in the table of [C-M, p.128]. 
Note that two non-conjugate parabolic subgroups have 
the same Richardson orbit $\mathcal{O}_{F_4(a_3)}$.  
By looking at the closure ordering of $F_4$ orbits [C, p.440], we 
see that the closure of each orbit contain a codimension 2 orbit.     

When $\g$ is of type $E_6$, there are exactly 6 single marked Dynkin 
diagrams. Four of them are already contained in the list of (ii). 
The Richardson orbits corresponding to other diagrams are 
$\mathcal{O}_{A_2}$ and $\mathcal{O}_{D_4(a_1)}$ in the list of 
$E_6$ nilpotent orbits in [C-M, p.129]. 
$\bar{\mathcal{O}}_{A_2}$ contains a codimension 2 orbit 
$\mathcal{O}_{3A_1}$. $\bar{\mathcal{O}}_{D_4(a_1)}$ contains 
a codimension 2 orbit $\mathcal{O}_{A_3 + A_1}$.       

When $\g$ is of type $E_7$, there are exactly 7 single marked Dynkin 
diagrams. Richardson orbits of the parabolic subgroups 
corresponding to them are $\mathcal{O}_{(3A_1)''}$, $\mathcal{O}_{A_2}$, 
$\mathcal{O}_{2A_2}$, $\mathcal{O}_{A_2 + 3A_1}$, $\mathcal{O}_{D_4(a_1)}$, 
$\mathcal{O}_{A_3 + A_2 + A_1}$ and $\mathcal{O}_{A_4 + A_2}$ in the 
table of [C-M, p.130-p.131]. 
By looking at the closure ordering of $E_7$ orbits [C, p.442], we 
see that the closure of each orbit contains a codimension 2 orbit. 

When $\g$ is of type $E_8$, there are exactly 8 single marked Dynkin 
diagrams. In the table of [C-M, 
p.132-p.134], Richardson orbits of the parabolic subgroups corresponding 
to them are $\mathcal{O}_{A_2}$, $\mathcal{O}_{2A_2}$, 
$\mathcal{O}_{D_4(a_1)}$, $\mathcal{O}_{D_4(a_1) + A_2}$, 
$\mathcal{O}_{A_4 + A_2}$, $\mathcal{O}_{A_4 + A_2 + A_1}$, 
$\mathcal{O}_{E_8(a_7)}$ and $\mathcal{O}_{A_6 + A_1}$. 
By looking at the closure ordering of $E_8$ orbits, we see that 
the closure of each orbit contains a codimension 2 orbit.    

STEP 2:   
Assume that $\g$ is classical. 
Let $f: \tilde{\mathcal{O}} \to \bar{\mathcal{O}}$ 
be the normalization map. By STEP 1 we may assume that 
$\bar{\mathcal{O}}$ contains a codimension 2 orbit $\mathcal{O}'$. 
In the classical case, by [K-P, 14], we see that $\tilde{\mathcal{O}}$ 
has actually singularities along $f^{-1}(\mathcal{O}')$.  
The Springer map $\mu$ is factorized as 
$$ T^*(G/P) \stackrel{\mu'}\to \tilde{\mathcal{O}} 
\stackrel{f}\to
 \mathcal{O}.$$ 
If $\mathrm{deg}(\mu) = 1$, then ${\mu}'$ is a birational maps 
of normal varieties. Then, by Zariski's main theorem, 
${\mu}'$ must have a positive dimensional fiber over a point 
of $f^{-1}(\mathcal{O}')$. This implies that $\mu$ is a 
divisorial birational map. 

Assume that $\g$ is of exceptional type. 
As explained above, the codimension 2 orbit 
$\mathcal{O}'$ of $\bar{\mathcal{O}}$ can 
be specified. It is enough to show that 
$\mathcal{O}'$ is $\mu$-relevant. By the previous 
proposition, we have to check that 
$V_{(x,1)}$ occurs in $\epsilon_{W(L)}^W$ for 
$x \in \mathcal{O}'$. In [Al], Alvis describes  
an irreducible decomposition of the induced 
representation $\mathrm{Ind}_{W(L)}^W(\rho)$ for  
any irreducible representation $\rho$ of $W(L)$. 
Hence, this can be done by using the tables of 
[Al] (see also the tables in [A-L], [B-L] and 
[C, 13.3]).    
Note that Spaltenstein [S] (cf. the footnote 
of p.68, [B-M]) has already checked that 
a special orbit is $\mu$-relevant by using 
these tables. Hence it is enough to check 
for non-special orbits $\mathcal{O}'$. One can 
find which orbits are non-special in the tables 
of [C-M, 8.4].

\begin{Exam}
({\bf Mukai flops}): Let $P$ and $P'$ be two parabolic subgroups of $G$ which 
correspond to dual marked Dynkin diagrams in the proposition 
above. Let $\mathcal{O}$ be the Richardson orbit of them. 
Then we have a diagram 
$$ T^*(G/P) \stackrel{\mu}\rightarrow \bar{\mathcal{O}} 
\stackrel{\mu'}\leftarrow T^*(G/P').$$ 
The birational maps $\mu$ and $\mu'$ are both small. 
Moreover, $T^*(G/P) - - \to T^*(G/P')$ is not an 
isomorphism. In fact, $T^*(G/P)$, $T^*(G/P')$ and 
$\bar{\mathcal{O}}$ all have $G$ actions, and $\mu$ and 
$\mu'$ are $G$-equivariant. If the birational map is 
an isomorphism, this would become a $G$-equivariant 
isomorphism. This implies that $G/P$ and $G/P'$ are 
isomorphic as $G$-varieties. In particular, $P$ and $P'$ 
are $G$-conjugate, which is absurd. Since the relative 
Picard numbers $\rho (T^*(G/P)/\bar{\mathcal{O}})$ and 
$\rho (T^*(G/P')/\bar{\mathcal{O}})$ equal $1$, we see 
that the diagram above is a {\em flop}. The diagram is 
called a {\em Mukai flop} of type $A_{n-1,k}$ (resp. $D_n$, 
$E_{6,I}$, $E_{6,II}$) according to the type of 
the corresponding marked Dynkin diagram.     
\end{Exam}        
 
\begin{Def} 
(i) Let $\mathcal{D}$ be a marked Dynkin diagram with 
exactly $l$ marked vertices. Choose $l-1$ marked vertices from 
them. Making the remained one vertex unmarked, we have a 
new marked Dynkin diagram $\bar{\mathcal{D}}$. This procedure 
is called a {\em contraction} of a marked Dynkin diagram. 
Next remove from $\mathcal{D}$ these $l-1$ vertices and 
edges touching these vertices. We then have a (non-connected) 
diagram; one of its connected component is a single marked 
Dynkin diagram. Assume that this single marked Dynkin diagram 
is one of those listed in Proposition 3.3. Replace this single 
marked Dynkin diagram by its dual and leave other components 
untouched.  Connecting again removed edges and vertices 
as before, we obtain a new marked Dynkin diagram $\mathcal{D'}$.  
Note that $\mathcal{D'}$ (resp. $\bar{\mathcal{D}}$) 
has exactly $l$ (resp. $l-1$) marked vertices. 
Now we say that $\mathcal{D'}$ is {\em adjacent} to 
$\mathcal{D}$ by means of $\bar{\mathcal{D}}$.   
\vspace{0.2cm}

(ii) Two marked Dynkin diagrams $\mathcal{D}$ and 
$\mathcal{D'}$ are called {\em equivalent} and 
are written as $\mathcal{D} \sim \mathcal{D'}$ if there 
is a finite chain of adjacent diagrams connecting 
$\mathcal{D}$ and $\mathcal{D'}$. 
\vspace{0.2cm}

(iii) Let $P$ be a parabolic subgroup of $G$ and 
let $\mathcal{D}_P$ be the corresponding marked 
Dynkin diagram. Two parabolic subgroups $P$ and $P'$ 
of $G$ are called {\em equivalent} and are written as 
$P \sim P'$ if $\mathcal{D}_P \sim \mathcal{D}_{P'}$.  
\end{Def}  

\begin{Exam}
Let us consider the marked Dynkin diagram 

\begin{picture}(300,20)(0,0) 
\put(0,0){$\mathcal{D}$:}
\put(30,0){\circle{5}}\put(35,0){\line(1,0){15}} 
\put(55,0){\circle*{5}}\put(55,-10){2}\put(60,0){\line(1,0){15}}
\put(80,0){\circle*{5}}\put(80,-10){3}\put(85,-3){$\Rightarrow$}
\put(100,0){\circle{5}}    
\end{picture}  
\vspace{0.3cm}

where vertices $2$ and $3$ are marked. 
We choose the vertex $3$. Making the remained 
one vertex (= the vertex $2$) unmarked, we have 
a marked Dynkin diagram 

\begin{picture}(300,20)(0,0) 
\put(0,0){$\bar{\mathcal{D}}$:}
\put(30,0){\circle{5}}\put(35,0){\line(1,0){15}} 
\put(55,0){\circle{5}}\put(55,-10){2}\put(60,0){\line(1,0){15}}
\put(80,0){\circle*{5}}\put(80,-10){3}\put(85,-3){$\Rightarrow$}
\put(100,0){\circle{5}}    
\end{picture}  
\vspace{0.3cm}

Now the following marked Dynkin diagram $\mathcal{D}'$ 
is adjacent to $\mathcal{D}$ by $\bar{\mathcal{D}}$. 

\begin{picture}(300,20)(0,0) 
\put(0,0){$\mathcal{D}'$:}
\put(30,0){\circle*{5}}\put(35,0){\line(1,0){15}} 
\put(55,0){\circle{5}}\put(55,-10){2}\put(60,0){\line(1,0){15}}
\put(80,0){\circle*{5}}\put(80,-10){3}\put(85,-3){$\Rightarrow$}
\put(100,0){\circle{5}}    
\end{picture}  
\vspace{0.3cm}
\end{Exam}

\section{Main Theorem} 
The following is our main theorem. For the notion of a 
relative ample cone and a relative movable cone, see 
[Ka], where some elementary roles of these cones in 
birational geometry are discussed. 

\begin{Thm}\label{main}
Let $\mathcal{O} \subset \g$ be a nilpotent 
orbit of a complex simple Lie algebra $\g$.  
Assume that 
its closure $\bar{\mathcal{O}}$ has a Springer 
resolution $\mu_{P_0}: T^*(G/P_0) \to \bar{\mathcal{O}}$. 
Then the following hold. \vspace{0.12cm}

(i) For a parabolic subgroup $P$ of $G$ such that 
$P \sim P_0$, $Y_P := T^*(G/P)$ gives a symplectic 
resolution of $\bar{\mathcal{O}}$. Conversely, any 
symplectic resolution is a Springer resolution of 
this form. \vspace{0.12cm}

(ii) The closure $\overline{\mathrm{Amp}}(Y_P/\bar{\mathcal{O}})$ 
of the relative ample cone is a simplicial polyhedral 
cone.  \vspace{0.12cm}

(iii) $\overline{\mathrm{Mov}}(Y_{P_0}/\overline{\mathcal{O}}) = 
\cup_{P \sim P_0} 
\overline{\mathrm{Amp}}(Y_P/\overline{\mathcal{O}})$, where 
$\overline{\mathrm{Mov}}(Y_{P_0}/\overline{\mathcal{O}})$ is 
the closure of the relative movable cone of $Y_{P_0}$ 
over $\overline{\mathcal{O}}$.  \vspace{0.12cm} 

(iv) A codimension $1$ face of $\overline{\mathrm{Amp}}(Y_P/ 
\overline{\mathcal{O}})$ corresponds to a small birational 
contraction map when it is a face of 
another ample cone, and corresponds to a divisorial 
contraction map when it is not a face of any other ample 
cone. \vspace{0.12cm}

(v) $\{Y_P\}_{P \sim P_0}$ are connected by Mukai flops 
of type $A$, $D$, $E_{6,I}$ and $E_{6,II}$.  
\end{Thm}

\begin{Rque}
For a classical complex Lie algebra, it is already known 
which nilpotent orbit closure has a Springer resolution 
(cf. [Na 2, Theorem 2.8]). When $\g$ is $G_2$, there are 
exactly 2 nilpotent orbits $\mathcal{O}_{G_2}$ and 
$\mathcal{O}_{G_2(a_1)}$ whose closures admit Springer 
resolutions. 
When $\g$ is $F_4$, such orbits are 
$\mathcal{O}_{A_2}$, $\mathcal{O}_{\tilde{A}_2}$, 
$\mathcal{O}_{F_4(a_3)}$, $\mathcal{O}_{B_3}$, 
$\mathcal{O}_{C_3}$, $\mathcal{O}_{F_4(a_2)}$, 
$\mathcal{O}_{F_4(a_1)}$ and $\mathcal{O}_{F_4}$. 
When $\g$ is $E_6$, such orbits are 
$\mathcal{O}_{2A_1}$, $\mathcal{O}_{A_2}$, $\mathcal{O}_{2A_2}$, 
$\mathcal{O}_{A_2 + 2A_1}$, $\mathcal{O}_{A_3}$, 
$\mathcal{O}_{D_4(a_1)}$, $\mathcal{O}_{A_4}$, 
$\mathcal{O}_{D_4}$, $\mathcal{O}_{A_4 + A_1}$, 
$\mathcal{O}_{D_5(a_1)}$, $\mathcal{O}_{E_6(a_3)}$, 
$\mathcal{O}_{D_5}$, $\mathcal{O}_{E_6(a_1)}$, and 
$\mathcal{O}_{E_6}$.    
\end{Rque}

The statement (ii) of Theorem 4.1 follows from 
the next Lemma.     

\begin{Lem} 
Let $G$ be a complex simple Lie group and let 
$P$ be a parabolic subgroup. Let $\hat{\mathcal{O}}$ 
be the Stein factorization of a Springer map 
$\mu: Y_P := T^*(G/P) \to \bar{\mathcal O}$. Then 
$\overline{Amp}(Y_P/\hat{\mathcal{O}})$ is a simplicial 
polyhedral cone. 
\end{Lem} 
 
{\em Proof}. Let $\mathcal{D}$ be the marked Dynkin 
diagram corresponding to $P$. Assume that $\mathcal{D}$ 
has $k$ marked vertices, say, $v_1$, ..., $v_k$. 
Then $b_2(G/P) = k$. Choose $l$ vertices 
$v_{i_1}, ..., v_{i_l}$, $1 \leq i_1 < ...< i_l \leq k$ 
and let $\mathcal{D}_{i_1, ..., i_l}$ be the marked 
Dynkin diagram such that exactly these $l$ vertices are 
marked and its underlying diagram is the same as $\mathcal{D}$. 
We denote by $X_{i_1, ..., i_l}$ the image of $Y_P \subset 
G/P \times \bar{\mathcal{O}}$ by the projection 
$$ G/P \times \bar{\mathcal{O}} \to 
G/P_{1_1, ..., i_l} \times \bar{\mathcal{O}}. $$ 
Let $$\nu_{i_1, ..., i_l}: Y_P \to X_{i_1, ..., i_l}$$ 
be the induced map. Then the Stein factorization of 
$\nu_{i_1, ..., i_l}$ is a birational contraction map, 
which corresponds to a codimension $k-l$ face of 
$\overline{Amp}(Y_P/\hat{\mathcal{O}})$. We shall denote 
by $F_{i_1, ..., i_l}$ this face. 
Then $\overline{\mathrm{Amp}}(Y_P/\hat
{\mathcal{O}})$ is a simplicial polyhedral cone 
generated by $F_1$, $F_2$, ..., and 
$F_k$. In fact, any $l$ dimensional face generated by 
$F_{i_1}$, ..., $F_{i_l}$ corresponds to the Stein 
factorization of $\nu_{i_1, ..., i_l}$, which is not 
an isomorphism.   
\vspace{0.2cm}

Next assume that two marked Dynkin diagrams 
$\mathcal{D}$ and $\mathcal{D'}$ are adjacent 
by means of $\bar{\mathcal{D}}$. We have three 
parabolic subgroups $P$, $P'$ and $\bar{P}$ of $G$ 
corresponding to $\mathcal{D}$, $\mathcal{D'}$ and 
$\bar{\mathcal{D}}$ respectively. One can assume that 
these subgroups contain the same Borel subgroup $B$ of 
$G$ and $\bar{P}$ contains both $P$ and $P'$. 
Let $\mu : T^*(G/P) \to \g$ and $\mu': T^*(G/P') \to 
\g$ be the Springer maps. 

\begin{Prop}\label{flop} 
(i) The Richardson orbits $\mathcal{O}$ of $P$ 
is the Richardson orbit of $P'$ 
\vspace{0.12cm}

(ii) Let $\nu$ be the composed map 
$$ T^*(G/P) \to G/P \times \bar{\mathcal{O}} 
\to G/\bar{P} \times \bar{\mathcal{O}} $$ and 
let $\nu'$ be the composed map 
$$ T^*(G/P') \to G/P' \times \bar{\mathcal{O}} 
\to G/\bar{P} \times \bar{\mathcal{O}}. $$  
Then $\mathrm{Im}(\nu) = \mathrm{Im}(\nu')$.  
\vspace{0.12cm}

(iii) If we put $X := \mathrm{Im}(\nu)$, then  
$$ T^*(G/P) \rightarrow X 
\leftarrow T^*(G/P') $$ 
is a locally trivial family of Mukai flops of type 
$A$, $D$, $E_{6,I}$ or $E_{6,II}$.  
In particular, $\nu$ and $\nu'$ are both small birational 
maps. If $\mathrm{deg}(\mu) = 1$, then $\mathrm{deg}(\mu') 
= 1$.       
\end{Prop}  

{\em Proof}. (i): Take a Levi decomposition 
$$ \bar{\mathfrak{p}} = l(\bar{\mathfrak{p}}) \oplus 
n(\bar{\mathfrak{p}}). $$ 
In the reductive Lie algebra $l(\bar{\mathfrak{p}})$, 
$\mathfrak{p} \cap l(\bar{\mathfrak{p}})$ and 
$\mathfrak{p}' \cap l(\bar{\mathfrak{p}})$ are parabolic 
subalgebras corresponding to dual marked Dynkin diagrams 
in Proposition 3.3. Hence they have conjugate Levi factors. 
On the other hand, we have 
$$ l(\mathfrak{p}) = l(\mathfrak{p} \cap l(\bar{\mathfrak{p}})),$$ 
and 
$$ l(\mathfrak{p}') = l(\mathfrak{p}' \cap l(\bar{\mathfrak{p}})).$$ 
Therefore, $l(\mathfrak{p})$ and $l(\mathfrak{p}')$ are 
conjugate. Since $\mathfrak{p}$ and $\mathfrak{p}'$ have 
conjugate Levi factors, their Richardson orbits coincide. 
\vspace{0.12cm}

(ii): Let $\mathcal{O}$ be the Richardson orbit of 
$\mathfrak{p}$ and $\mathfrak{p}'$. 
Springer maps $\mu : T^*(G/P) \to \bar{\mathcal{O}}$ 
and $\mu' : T^*(G/P') \to \bar{\mathcal{O}}$ are both 
$G$-equivariant with respect to natural $G$-actions. 
Then $U := \mu^{-1}(\mathcal{O})$ and 
$U' := ({\mu}')^{-1}(\mathcal{O})$ are open dense 
orbits of $T^*(G/P)$ and $T^*(G/P')$ respectively. 
Since $\nu$ and ${\nu}'$ are proper maps, 
$\mathrm{Im}(\nu) = \overline{\nu(U)}$ and 
$\mathrm{Im}({\nu}') = \overline{\nu'(U')}$.  
In the following we shall prove that $\nu(U) = 
{\nu}'(U')$.

(ii-1): 
We regard $T^*(G/P)$ (resp. $T^*(G/P')$) as 
a closed subvariety of $G/P \times \bar{\mathcal{O}}$ 
(resp. $G/P' \times \bar{\mathcal{O}}$). 
By replacing $P'$ by a suitable conjugate in 
$\bar{P}$, we may assume that there exists an element 
$x \in \mathcal{O}$ such that $([P], x) \in U$ 
and $([P'], x) \in U'$. In fact, for a Levi decomposition 
$$ \bar{\mathfrak{p}} = l(\bar{\mathfrak{p}}) \oplus 
n(\bar{\mathfrak{p}}),$$ we have a direct sum decomposition 
$$ n(\mathfrak{p}) = n(\mathfrak{p} \cap l(\bar{\mathfrak{p}})) 
\oplus n(\bar{\mathfrak{p}}).$$ 
Let $p_1: n(\mathfrak{p}) \to n(\mathfrak{p} \cap 
l(\bar{\mathfrak{p}}))$ be the 1-st projection. 
Let $\mathcal{O}' \subset l(\bar{\mathfrak{p}})$ be the 
Richardson orbit of the parabolic 
subalgebra $\mathfrak{p} \cap l(\bar{\mathfrak{p}})$ 
of $l(\bar{\mathfrak{p}})$. 
Since $p_1^{-1}(n(\mathfrak{p}) \cap \mathcal{O}')$ 
and $n(\mathfrak{p}) \cap \mathcal{O}$ are both Zariski 
open subsets of $n(\mathfrak{p})$, we can take an element 
$$ x \in p_1^{-1}(n(\mathfrak{p}) \cap \mathcal{O}') \cap 
(n(\mathfrak{p}) \cap \mathcal{O}). $$ 
Since $x \in n(\mathfrak{p}) \cap \mathcal{O}$, we have 
$([P],x) \in U$. Decompose $x = x_1 + x_2$ according to the 
direct sum decomposition. Then $x_1 \in \mathcal{O}'$. 
The orbit $\mathcal{O}'$ is also the Richardson orbit of 
$\mathfrak{p}' \cap l(\bar{\mathfrak{p}})$. Therefore, 
for some $g \in L(\bar{P})$ (the Levi factor of $\bar{P}$ 
corresponding to $l(\bar{P})$), 
$$x_1 \in n(Ad_g(\mathfrak{p}' 
\cap l(\bar{\mathfrak{p}}))).$$ 
The Levi decomposition of $\bar{\mathfrak{p}}$ induces 
a direct sum decomposition 
$$ n(Ad_g(\mathfrak{p}')) = n(Ad_g(\mathfrak{p}') 
\cap l(\bar{\mathfrak{p}})) \oplus n(\bar{\mathfrak{p}}).$$ 
Note that $Ad_g(\mathfrak{p}') \cap l(\bar{\mathfrak{p}}) 
= Ad_g(\mathfrak{p}' \cap l(\bar{\mathfrak{p}}))$. 
Hence we see that $x_1 + x_2 \in n(Ad_g(\mathfrak{p}'))$. 
Now, for $Ad_g(P') \subset \bar{P}$, we have 
$([Ad_g(P')], x) \in U'$. 

(ii-2): Any element of $U$ can be written as 
$([gP], Ad_g(x))$ for some $g \in G$. Then 
$$\nu([gP], Ad_g(x)) = ([g\bar{P}], Ad_g(x)).$$ 
For the same $g \in G$, we have 
$([gP'], Ad_g(x)) \in U'$ and 
$$\nu'([gP'], Ad_g(x)) = ([g\bar{P}], Ad_g(x)).$$ 
Therefore, $\nu(U) \subset \nu'(U')$. By the same 
argument, we also have $\nu'(U') \subset \nu(U)$.   
\vspace{0.12cm}

(iii): For $g \in G$, $Ad_g(n(\bar{\mathfrak{p}}))$ 
is the nil-radical of $Ad_g(\bar{\mathfrak{p}})$. 
Since $Ad_g(\bar{\mathfrak{p}})$ depends only on 
the class $[g] \in G/\bar{P}$, $Ad_g(n(\bar{\mathfrak{p}}))$ 
also depends on the class $[g] \in G/\bar{P}$. 
We denote by $Ad_g(l(\bar{\mathfrak{p}}))$ the quotient 
of $Ad_g(\bar{\mathfrak{p}})$ by its nil-radical 
$Ad_g(n(\bar{\mathfrak{p}}))$. 
Let us consider the vector bundle over $G/\bar{P}$ 
$$ \cup_{[g] \in G/\bar{P}}Ad_g(\bar{\mathfrak{p}}) 
\to G/\bar{P}.$$  
Let $\mathcal{L}$ be its quotient bundle whose fiber 
over $[g] \in G/\bar{P}$ is $Ad_g(l(\bar{\mathfrak{p}}))$. 
We call $\mathcal{L}$ the Levi bundle. 
Let $\mathcal{O}'$ be the Richardson orbit of the parabolic 
subalgebra $\mathfrak{p} \cap l(\bar{\mathfrak{p}})$ of 
$l(\bar{\mathfrak{p}})$. Note that $\mathcal{O}'$ is also 
the Richardson orbit of $\mathfrak{p}' \cap 
l(\bar{\mathfrak{p}})$. 
In $\mathcal{L}$, we consider the fiber bundle  
$$ W := \cup_{[g] \in G/\bar{P}}Ad_g(\bar{\mathcal{O}'}) $$ 
whose fiber over 
$[g] \in G/\bar{P}$ is $Ad_g(\bar{\mathcal{O}'})$.  
Put $X := \mathrm{Im}(\nu)$. Define a map  
$$ f : X \to W $$ 
as $f([g], x) := ([g], x_1)$, where 
$x_1$ is the first factor of $x$ under the 
direct sum decomposition 
$$ Ad_g(\bar{\mathfrak{p}}) = Ad_g(l(\bar{\mathfrak{p}})) 
\oplus n(Ad_g(\bar{\mathfrak{p}})).$$ 
Note that $x_1 \in Ad_g(\bar{\mathcal{O}'})$. 
In fact, in the direct sum decomposition, 
we have 
$$ n(Ad_g(\mathfrak{p})) = n(Ad_g(\mathfrak{p}) \cap 
Ad_g(l(\bar{\mathfrak{p}}))) \oplus n(Ad_g(\bar{\mathfrak{p}})).$$ 
Therefore 
$$ x_1 \in n(Ad_g(\mathfrak{p}) \cap 
Ad_g(l(\bar{\mathfrak{p}}))) \subset Ad_g(\bar{\mathcal{O}'}).$$ 
Since $W \to G/\bar{P}$ is an $\bar{\mathcal{O}'}$ bundle, we 
have a family of Mukai flops parametrized by $G/\bar{P}$: 
$$ Y \rightarrow W \leftarrow Y'.$$ 
By pulling back this diagram by $f: X \to W$, we 
have the diagram 
$$ T^*(G/P) \rightarrow X \leftarrow T^*(G/P').$$ 
\vspace{0.2cm}

Let $\mathcal{D}$ be a marked Dynkin diagram 
and let $\bar{\mathcal{D}}$ be the diagram obtained 
from $\mathcal{D}$ by a contraction. 
Let $P$ and $\bar{P}$ be parabolic subgroups of $G$ 
corresponding to $\mathcal{D}$ and $\bar{\mathcal{D}}$ 
respectively. One can assume that $\bar{P}$ contains 
$P$. Let $\mathcal{O}$ be the Richardson orbit 
of $P$ and let $\nu$ be the compoed map 
$$ T^*(G/P) \to G/P \times \bar{\mathcal{O}} 
\to G/\bar{P} \times \bar{\mathcal{O}}. $$ 
We put $X := \mathrm{Im}(\nu)$. As above, $\mu: 
T^*(G/P) \to \bar{\mathcal{O}}$ is the Springer map. 

\begin{Prop}\label{dvisorial}
Let $\g$ be a complex simple Lie algebra. 
Assume that no marked Dynkin diagram is 
adjacent to $\mathcal{D}$ by means of 
$\bar{\mathcal{D}}$. If $\mathrm{deg}(\mu) 
= 1$, then $\nu: T^*(G/P) \to X$ is a 
divisorial birational contraction map. 
\end{Prop}  

{\em Proof}. As in the proof of Proposition 4.4, (iii), 
we construct an $\bar{\mathcal{O}'}$ bundle 
$W$ over $G/\bar{P}$ and define a map 
$f: X \to W.$  There is a family of 
Springer maps 
$$ Y \stackrel{\sigma}\to W \to G/\bar{P}. $$ 
By pulling back $Y \stackrel{\sigma}\to W$ 
by $f: X \to W$, we have the $\nu: T^*(G/P) 
\to X$. Since $\mathrm{deg}\mu = 1$, $\nu$ is 
a birational map. Hence $\sigma$ should be a 
birational map. Hence $\sigma : Y \to W$ is a 
family of Springer resolutions. By the assumption, 
there are no marked Dynkin diagrams adjacent to 
$\mathcal{D}$ by means of $\bar{\mathcal{D}}$. 
Now Proposition 3.3 shows 
that the Springer resolution is divisorial. 
Therefore, $\nu$ is also divisorial.  
\vspace{0.2cm}

Now let us prove Theorem 4.1. By Proposition 
4.4, (iii), $Y_P := T^*(G/P)$ all give symplectic 
resolutions of $\bar{\mathcal{O}}$ for $P \sim P_0$. 
Hence the first statement of (i) has been proved. 
Moreover, $\{Y_P\}$ are connected by Mukai flops, 
which is nothing but (v). Let us consider 
$\cup_{P \sim P_0}\overline{Amp}(Y_P/\bar{\mathcal{O}})$ 
in $N^1(Y_{P_0}/\bar{\mathcal{O}})$. Then  (iv) 
follows from Proposition 4.4, (iii) and Proposition 4.5.   
For an $\bar{\mathcal{O}}$-movable divisor $D$ on $Y_{P_0}$, 
a $K_{Y_{P_0}} + D$-extremal contraction is a small birational 
map. Therefore, the corresponding codimension $1$ face of 
$\overline{\mathrm{Amp}}(Y_{P_0}/\bar{\mathcal{O}})$ becomes 
a codimension 1 face of another $\overline{\mathrm{Amp}}(Y_P/ 
\bar{\mathcal{O}})$. For this small birational map, there 
exists a flop. Replace $D$ by its proper transform and 
continue the same. We shall prove that this procedure 
ends in finite times. Suppose to the contrary. Since the flops 
occur between finite number of varieties $\{Y_P\}$, 
a variety, say $Y_{P_1}$, appears at least twice in the sequence 
of flops: 
$$ Y_{P_1} --\to Y_{P_2} --\to ... --\to Y_{P_1}. $$ 
For the first flop 
$$ Y_{P_1} \stackrel{\nu_1}\rightarrow 
X_1 \leftarrow Y_{P_2},$$ 
take a discrete valation $v$ of the function field 
$K(Y_{P_1})$  in such a way that its center is contained 
in the exceptional locus $\mathrm{Exc}(\nu_1)$ of $\nu_1$. 
Let $D_i \subset Y_{P_i}$ be the proper transforms 
of $D$. Then we have inequalities for {\em discrepancies}
(cf. [KMM], Proposition 5-1-11): 
$$ a(v, D_1) < a(v, D_2) \leq .... \leq a(v, D_1). $$ 
Here the first inequality is a strict one since 
the center of $v$ is contained in $\mathrm{Exc}(\nu_1)$. 
This is absurd. Hence the procedure ends in finite times, 
which implies that $D \in \overline{\mathrm{Amp}}(Y_P/
\bar{\mathcal{O}})$ for some $P$. Therefore, (iii) has been 
proved. The second statement of (i) immediately follows 
from (iii).    
  
\begin{Exam}
([Na 2, Example 4.6]): Assume that $\g = \mathfrak{sl}(6)$. 
The marked Dynkin diagram $\mathcal{D}$ 

\begin{picture}(300,20)
\put(30,0){\circle*{5}}\put(35,0){\line(1,0){20}}
\put(60,0){\circle{5}}\put(65,0){\line(1,0){20}}
\put(90,0){\circle*{5}}\put(95,0){\line(1,0){20}}
\put(120,0){\circle{5}}\put(125,0){\line(1,0){20}}
\put(150,0){\circle{5}} 
\end{picture} 
\vspace{0.4cm}

gives a parabolic subgroup $P_{1,2,3} \subset SL(6)$ 
of flag type $(1,2,3)$. We put $Y_{1,2,3} := 
T^*(G/P_{1,2,3})$. There are 5 other marked Dynkin 
diagrams which are equivalent to $\mathcal{D}$: 

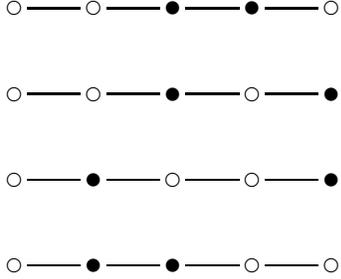
\begin{picture}(300,20)
\put(30,0){\circle*{5}}\put(35,0){\line(1,0){20}}
\put(60,0){\circle{5}}\put(65,0){\line(1,0){20}}
\put(90,0){\circle{5}}\put(95,0){\line(1,0){20}}
\put(120,0){\circle*{5}}\put(125,0){\line(1,0){20}}
\put(150,0){\circle{5}} 
\end{picture} 
\vspace{0.4cm}

\begin{picture}(300,20)
\put(30,0){\circle{5}}\put(35,0){\line(1,0){20}}
\put(60,0){\circle{5}}\put(65,0){\line(1,0){20}}
\put(90,0){\circle*{5}}\put(95,0){\line(1,0){20}}
\put(120,0){\circle*{5}}\put(125,0){\line(1,0){20}}
\put(150,0){\circle{5}} 
\end{picture} 
\vspace{0.4cm}

\begin{picture}(300,20)
\put(30,0){\circle{5}}\put(35,0){\line(1,0){20}}
\put(60,0){\circle{5}}\put(65,0){\line(1,0){20}}
\put(90,0){\circle*{5}}\put(95,0){\line(1,0){20}}
\put(120,0){\circle{5}}\put(125,0){\line(1,0){20}}
\put(150,0){\circle*{5}} 
\end{picture} 
\vspace{0.4cm}

\begin{picture}(300,20)
\put(30,0){\circle{5}}\put(35,0){\line(1,0){20}}
\put(60,0){\circle*{5}}\put(65,0){\line(1,0){20}}
\put(90,0){\circle{5}}\put(95,0){\line(1,0){20}}
\put(120,0){\circle{5}}\put(125,0){\line(1,0){20}}
\put(150,0){\circle*{5}} 
\end{picture} 
\vspace{0.4cm}

\begin{picture}(300,20)
\put(30,0){\circle{5}}\put(35,0){\line(1,0){20}}
\put(60,0){\circle*{5}}\put(65,0){\line(1,0){20}}
\put(90,0){\circle*{5}}\put(95,0){\line(1,0){20}}
\put(120,0){\circle{5}}\put(125,0){\line(1,0){20}}
\put(150,0){\circle{5}} 
\end{picture} 
\vspace{0.4cm}

Five parabolic subgroups $P_{1,3,2}$, $P_{3,1,2}$, 
$P_{3,2,1}$, $P_{2,3,1}$, $P_{2,1,3}$ correspond 
to the marked Dynkin diagrams above respectively. 
We put $Y_{i,j,k} := T^*(SL(6)/P_{i,j,k})$.  
Let $\mathcal{O}$ be the Richardson orbit of these 
parabolic subgroups. Then $\overline{\mathrm{Mov}}(Y_{1,2,3}/
\bar{\mathcal{O}}) \cong \mathbf{R}^2$, which is divided 
into six chambers by the ample cones of $Y_{i,j,k}$ in the 
following way:

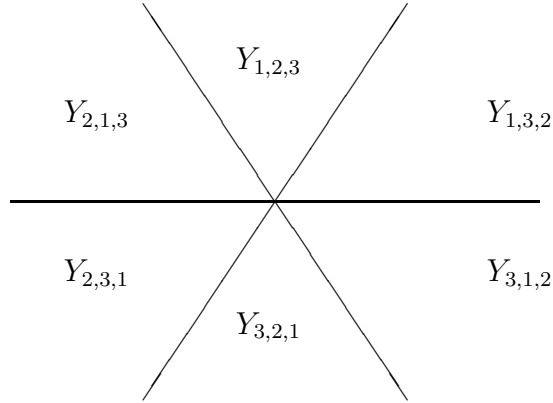
\begin{picture}(300,100)(0,0) 
\put(150,0){\line(2,3){50}}\put(150,0){\line(1,0){100}}
\put(150,0){\line(2,-3){50}}\put(150,0){\line(-2,-3){50}}
\put(150,0){\line(-1,0){100}}\put(150,0){\line(-2,3){50}}
\put(135,50){$Y_{1,2,3}$}\put(230,30){$Y_{1,3,2}$}
\put(230,-30){$Y_{3,1,2}$}\put(135,-50){$Y_{3,2,1}$}
\put(70,-30){$Y_{2,3,1}$}\put(70,30){$Y_{2,1,3}$}
\end{picture} 
\vspace{3.0cm}
\end{Exam}

\begin{Exam}
([Na 2, Example 4.7]): 
Assume that $\g = \mathfrak{so}(10)$. 
The marked Dynkin diagram 

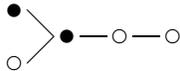
\begin{picture}(300,20)(0,0) 
\put(30,10){\circle*{5}}\put(30,-10){\circle{5}}
\put(35,10){\line(1,-1){10}}\put(35,-10){\line(1,1){10}}
\put(50,0){\circle*{5}}\put(55,0){\line(1,0){10}}
\put(70,0){\circle{5}}\put(75,0){\line(1,0){10}}
\put(90,0){\circle{5}}
\end{picture} 
\vspace{0.4cm}

gives a parabolic subgroup $P^+_{3,2,2,3}$ of flag 
type $(3,2,2,3)$. There are three marked Dynkin diagrams 
equivalent to this marked diagram: 

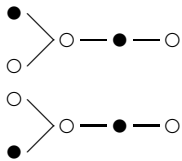
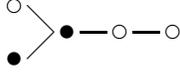
\begin{picture}(300,20)(0,0) 
\put(30,10){\circle*{5}}\put(30,-10){\circle{5}}
\put(35,10){\line(1,-1){10}}\put(35,-10){\line(1,1){10}}
\put(50,0){\circle{5}}\put(55,0){\line(1,0){10}}
\put(70,0){\circle*{5}}\put(75,0){\line(1,0){10}}
\put(90,0){\circle{5}}
\end{picture} 
\vspace{0.4cm}

\begin{picture}(300,20)(0,0) 
\put(30,10){\circle{5}}\put(30,-10){\circle*{5}}
\put(35,10){\line(1,-1){10}}\put(35,-10){\line(1,1){10}}
\put(50,0){\circle{5}}\put(55,0){\line(1,0){10}}
\put(70,0){\circle*{5}}\put(75,0){\line(1,0){10}}
\put(90,0){\circle{5}}
\end{picture} 
\vspace{0.4cm}

\begin{picture}(300,20)(0,0) 
\put(30,10){\circle{5}}\put(30,-10){\circle*{5}}
\put(35,10){\line(1,-1){10}}\put(35,-10){\line(1,1){10}}
\put(50,0){\circle*{5}}\put(55,0){\line(1,0){10}}
\put(70,0){\circle{5}}\put(75,0){\line(1,0){10}}
\put(90,0){\circle{5}}
\end{picture} 
\vspace{0.4cm}

Three parabolic subgroups $P^+_{2,3,3,2}$, $P^{-}_{2,3,3,2}$, 
$P^{-1}_{3,2,2,3}$ correspond 
to these marked Dynkin diagrams respectively. 
Note that there are exactly two conjugacy classes of 
parabolic subgroups with the same flag type (cf. 
Example 2.1). 
We put $Y^+_{i,j} := T^*(SO(10)/P^+_{i,j,j,i})$ 
and put $Y^{-}_{i,j} := T^*(SO(10)/P^{-}_{i,j,j,i}).$ 
Let $\mathcal{O}$ be the Richardson orbit of these 
parabolic subgroups. Then $\overline{\mathrm{Mov}}(Y^+_{3,2}/
\bar{\mathcal{O}})$ 
is divided into four chambers by the ample cones of $Y^{+}_{3,2}$, 
$Y^{+}_{2,3}$, $Y^{-}_{2,3}$, $Y^{-}_{3,2}$ in the 
following way:

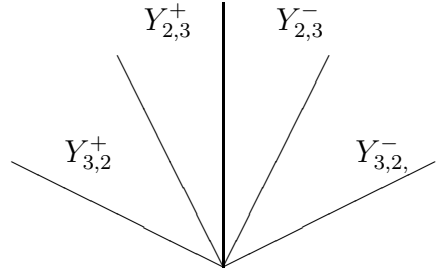
\begin{picture}(300,150)(0,0) 
\put(150,0){\line(0,1){100}}\put(150,0){\line(2,1){80}}
\put(150,0){\line(-2,1){80}}\put(150,0){\line(1,2){40}}
\put(150,0){\line(-1,2){40}}\put(170,90){$Y^{-}_{2,3}$}
\put(120,90){$Y^+_{2,3}$}\put(200,40){$Y^{-}_{3,2,}$}
\put(90,40){$Y^+_{3,2}$}
\end{picture} 
\vspace{0.4cm}
\end{Exam}

\begin{Exam}
Assume that $\g$ is of type $E_6$. 
Consider the nilpotent orbit $\mathcal{O} := 
\mathcal{O}_{A_3}$ 
(cf. [C-M], p.129). 
This is the unique orbit with dimension 52. 
By a dimension count, we see that $\mathcal{O}$ 
is the Richardson orbit of the parabolic subgroup 
$P_1 \subset G$ 
associated with the marked Dynkin diagram 

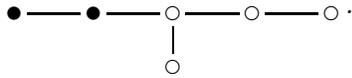
\begin{picture}(300,20)
\put(30,0){\circle*{5}}\put(35,0){\line(1,0){20}}
\put(60,0){\circle*{5}}\put(65,0){\line(1,0){20}}
\put(90,0){\circle{5}}\put(90,-5){\line(0,-1){10}}
\put(90,-20){\circle{5}}\put(95,0){\line(1,0){20}}
\put(120,0){\circle{5}}\put(125,0){\line(1,0){20}}
\put(150,0){\circle{5}.} 
\end{picture} 
\vspace{0.7cm}

Since $\pi_1(\mathcal{O}) = 1$ ([C-M], p.129], 
the Springer map $\nu_1: T^*(G/P_1) \to \bar{\mathcal{O}}$ 
has degree 1. The following marked Dynkin diagrams are 
equivalent to the diagram above: 

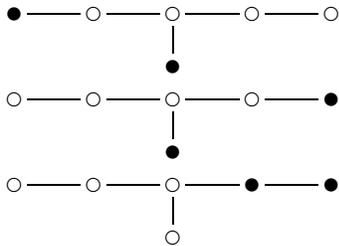
\begin{picture}(300,20)
\put(30,0){\circle*{5}}\put(35,0){\line(1,0){20}}
\put(60,0){\circle{5}}\put(65,0){\line(1,0){20}}
\put(90,0){\circle{5}}\put(90,-5){\line(0,-1){10}}
\put(90,-20){\circle*{5}}\put(95,0){\line(1,0){20}}
\put(120,0){\circle{5}}\put(125,0){\line(1,0){20}}
\put(150,0){\circle{5}} 
\end{picture} 
\vspace{0.4cm}

\begin{picture}(300,20)
\put(30,0){\circle{5}}\put(35,0){\line(1,0){20}}
\put(60,0){\circle{5}}\put(65,0){\line(1,0){20}}
\put(90,0){\circle{5}}\put(90,-5){\line(0,-1){10}}
\put(90,-20){\circle*{5}}\put(95,0){\line(1,0){20}}
\put(120,0){\circle{5}}\put(125,0){\line(1,0){20}}
\put(150,0){\circle*{5}} 
\end{picture} 
\vspace{0.4cm}

\begin{picture}(300,20)
\put(30,0){\circle{5}}\put(35,0){\line(1,0){20}}
\put(60,0){\circle{5}}\put(65,0){\line(1,0){20}}
\put(90,0){\circle{5}}\put(90,-5){\line(0,-1){10}}
\put(90,-20){\circle{5}}\put(95,0){\line(1,0){20}}
\put(120,0){\circle*{5}}\put(125,0){\line(1,0){20}}
\put(150,0){\circle*{5}.} 
\end{picture} 
\vspace{0.7cm}

Denote by $P_2$, $P_3$, $P_4$ respectively 
the parabolic subgroups 
corresponding to the diagrams above. 
We put $Y_i := T^*(G/P_i)$ for $i = 1,2,3,4$. 
Then $\overline{\mathrm{Mov}}(Y_1/\bar{\mathcal{O}})$ 
is divided into four chambers by the ample cones of 
$Y_i$: 

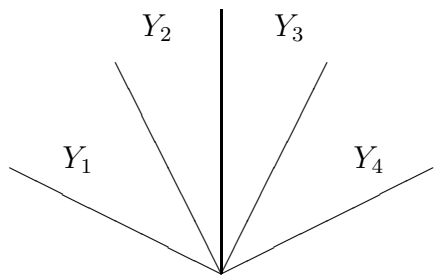
\begin{picture}(300,150)(0,0) 
\put(150,0){\line(0,1){100}}\put(150,0){\line(2,1){80}}
\put(150,0){\line(-2,1){80}}\put(150,0){\line(1,2){40}}
\put(150,0){\line(-1,2){40}}\put(170,90){$Y_3$}
\put(120,90){$Y_2$}\put(200,40){$Y_4$}
\put(90,40){$Y_1$}
\end{picture} 
\vspace{0.4cm}

$Y_1$ and $Y_2$ are connected by a Mukai flop of type 
$D_5$ (cf. Proposition 4.4, (iii)). 
$Y_2$ and $Y_3$ are connected by 
a Mukai flop of type $A_{5,1}$ (for the notation, see 
Example 2.5). $Y_3$ and $Y_4$ are connected by a Mukai 
flop of type $D_5$. 
\end{Exam}

\quad \\
\quad\\

Yoshinori Namikawa \\
Departement of Mathematics, 
Graduate School of Science, Kyoto University,
Kita-shirakawa Oiwake-cho, Kyoto, 606-8502, JAPAN \\
namikawa@kusm.kyoto-u.ac.jp
      
\end{document}